\documentclass{article}
\author{Johan Ernest Mebius\thanks{Delft University of Technology, Faculty of
Electrical Engineering, Mathematics and Computer Science, P.O.Box 5031, NL -- 2600 GA Delft,
The Netherlands, Phone +31.15.2783072, Fax +31.15.2787022}}
\title{DERIVATION OF THE EULER--RODRIGUES FORMULA\\ FOR THREE-DIMENSIONAL ROTATIONS\\
FROM THE GENERAL FORMULA\\ FOR FOUR-DIMENSIONAL ROTATIONS}
\date{January 2007}
\begin{document}
\maketitle

\begin{abstract}
\noindent
The general 4D rotation matrix is specialised to the general 3D rotation matrix by equating its $a_{00}$ element to $1$. Its associate matrix of products of the left--hand and right--hand quaternion components is specialised correspondingly. Inequalities involving the angles through which the coordinate axes in 3D space are displaced are used to prove that the left--hand and the right--hand quaternions are each other's inverses, thus proving the {\sc Euler--Rodrigues} formula.
\\
A general procedure to determine the {\sc Euler} parameters of a given 3D rotation matrix is sketched.
\\
By putting $a_{00} = -1$ instead of $+1$ in the general 4D rotation matrix one proves the counterpart of the {\sc Euler--Rodrigues} formula for 3D rotoreflections.
\\
\\
{\bf Keywords:} Euler--Rodrigues formula, Euler parameters, quaternions, four--dimensional rotations, three--dimensional rotations, rotoreflections
\end{abstract}

\section*{Notations}

\begin{tabbing}
WWWWWWWWWWJ\=DUMMY     \kill
$a, b, c, d$           \>In Section \ref{Section4}: {\sc Euler} parameters of a 3D rotation\\
$A$                    \>Arbitrary 4D rotation matrix\\
$M$                    \>Associate of matrix $A$ (defined in [MEBI~2005])\\
$M_L$, $M_R$           \>Matrices representing left-- and right--multiplication by a unit quaternion, respectively\\
$L = a + bi + cj + dk$ \>Unit quaternion appearing as a leftmost factor in quaternion multiplication\\
$P = u + xi + yj + zk$ \>Arbitrary 4D point represented as a quaternion\\
$R = p + qi + rj + sk$ \>Unit quaternion appearing as a rightmost factor in quaternion multiplication\\
$OXYZ$                 \>3D Cartesian coordinate system\\
$\rho$                 \>Distance of point in 3D space from $Z$ axis 
\end{tabbing}

\section{The quaternion representation theorem for 4D rotations}
\label{Section1}

{\sc Theorem:} Each 4D rotation matrix can be decomposed in two
ways into a matrix representing left--multiplication by a unit quaternion
and a matrix representing right--multiplication by a unit quaternion.
These decompositions differ only in the signs of the component matrices.
\\
In quaternion representation: $P' = LPR$; in matrix representation: $P' = M_LM_RP$.
\\
\\
{\sc Outline of proof:} Let $M_L$, $M_R$ (Eq. \ref{Label1}) be matrices representing
left-- and right--multiplication by a unit quaternion, respectively. Then
their product $A = M_LM_R$ (Eq. \ref{Label2}) is a 4D rotation matrix.
\\
Matrix $M_L$ is determined by four reals $a, b, c, d$ satisfying the relation
$a^2 + b^2 + c^2 + d^2 = 1$.
Likewise, matrix $M_R$ is determined by four reals $p, q, r, s$ satisfying the
relation $p^2 + q^2 + r^2 + s^2 = 1$.
\\
The 16 products $ap,\ aq,\ ar,\ as$, $\ldots$ , $dp,\ dq,\ dr,\ ds$ are arranged
into a matrix $M$ (Eq. \ref{Label4}), which has rank one and norm unity (when considered 
as a 16D vector), and is easily expressed in the elements of $A$. In [MEBI~2005] 
matrix $M$ is denoted as the {\em associate matrix} of $A$.
\\
Conversely, given an arbitrary 4D rotation matrix $A$ (Eq. \ref{Label3}), one calculates its
associate matrix $M$ according to Eq. \ref{Label5} in the hope that it is a matrix of products $ap,\ aq,\ ar,\ as$,
$\ldots$ , $dp,\ dq,\ dr,\ ds$ which are not all zero. This hope is vindicated
by proving that $M$ has rank one whenever $A$ is a 4D rotation matrix. For this proof one needs the general theorem that complementary subdeterminants of a rotation matrix of any number of dimensions are equal.
\\
The proof is completed by observing that the sum of the squares of the elements
of $M$ is unity, and concluding that two pairs of quadruples of reals
$a,\ b,\ c,\ d$;\ \ $p,\ q,\ r,\ s$ exist satisfying $a^2 + b^2 + c^2 + d^2 = 1$,\ \
$p^2 + q^2 + r^2 + s^2 = 1$ and differing only in sign.
\\
A complete proof is given in [MEBI~2005].
\\
\\
\begin{equation}
\label{Label1}
M_L=\left [\begin {array}{rrrr} a&-b&-c&-d\\\noalign{\medskip}b&a&-d&c\\\noalign{\medskip}c&d&a&-b\\\noalign{\medskip}d&-
c&b&a\end {array}\right ],
\ \ M_R=\left [\begin {array}{rrrr} p&-q&-r&-s\\\noalign{\medskip}q&p&s&-r\\\noalign{\medskip}r&-s&p&q\\\noalign{\medskip}s&r
&-q&p\end {array}\right ]
\end{equation}
\begin{equation}
\label{Label2}
M_LM_R =
\left[
\begin {array}{rrrr} 
                  ap-bq-cr-ds&-aq-bp+cs-dr&-ar-bs-cp+dq&-as+br-cq-dp\\
\noalign{\medskip}bp+aq-dr+cs&-bq+ap+ds+cr&-br+as-dp-cq&-bs-ar-dq+cp\\
\noalign{\medskip}cp+dq+ar-bs&-cq+dp-as-br&-cr+ds+ap+bq&-cs-dr+aq-bp\\
\noalign{\medskip}dp-cq+br+as&-dq-cp-bs+ar&-dr-cs+bp-aq&-ds+cr+bq+ap
\end {array}
\right]
\end{equation}
\begin{equation}
\label{Label3}
A =
\left[
\begin{array}{rrrr}
                    a_{00}  &   a_{01}  &   a_{02}  &   a_{03}  \\
\noalign{\medskip}  a_{10}  &   a_{11}  &   a_{12}  &   a_{13}  \\
\noalign{\medskip}  a_{20}  &   a_{21}  &   a_{22}  &   a_{23}  \\
\noalign{\medskip}  a_{30}  &   a_{31}  &   a_{32}  &   a_{33}
\end{array}
\right]
\end{equation}
\begin{equation}
\label{Label4}
M =
\left[
\begin{array}{rrrr}
                    ap  &   aq  &   ar  &   as  \\
\noalign{\medskip}  bp  &   bq  &   br  &   bs  \\
\noalign{\medskip}  cp  &   cq  &   cr  &   cs  \\
\noalign{\medskip}  dp  &   dq  &   dr  &   ds
\end{array}
\right]
=
\left[
\begin{array}{rrrr}
                    a  \\
\noalign{\medskip}  b  \\
\noalign{\medskip}  c  \\
\noalign{\medskip}  d
\end{array}
\right]
\begin{array}{r}
\left[
\begin{array}{rrrr}
p  &   q  &   r  &   s
\end{array}
\right]
\\
\begin{array}{rrrr}
\   & \   & \   & \   \\
\   & \   & \   & \   \\
\   & \   & \   & \   \\
\   & \   & \   & \
\end{array}
\end{array}
\end{equation}
\begin{equation}
\label{Label5}
= \frac{1}{4}
\left[
\begin{array}{rrrr}
                   a_{00}+a_{11}+a_{22}+a_{33} & +a_{10}-a_{01}-a_{32}+a_{23} & +a_{20}+a_{31}-a_{02}-a_{13} & +a_{30}-a_{21}+a_{12}-a_{03} \\
\noalign{\medskip} a_{10}-a_{01}+a_{32}-a_{23} & -a_{00}-a_{11}+a_{22}+a_{33} & +a_{30}-a_{21}-a_{12}+a_{03} & -a_{20}-a_{31}-a_{02}-a_{13} \\
\noalign{\medskip} a_{20}-a_{31}-a_{02}+a_{13} & -a_{30}-a_{21}-a_{12}-a_{03} & -a_{00}+a_{11}-a_{22}+a_{33} & +a_{10}+a_{01}-a_{32}-a_{23} \\
\noalign{\medskip} a_{30}+a_{21}-a_{12}-a_{03} & +a_{20}-a_{31}+a_{02}-a_{13} & -a_{10}-a_{01}-a_{32}-a_{23} & -a_{00}+a_{11}+a_{22}-a_{33}
\end{array}
\right]
\end{equation}
\\

\section{3D Rotations}
\label{Section2}

Any 3D isometry with a fixed point $O$ has an invariant plane and an invariant axis which intersect at right angles at $O$. 
\\
An orientation--preserving isometry with fixed point $O$ is a 3D rotation through a certain angle $\alpha$ around its invariant axis in its invariant plane. (Theorem of {\sc Euler})
\\
An orientation--reversing isometry with fixed point $O$ is a 3D rotation through a certain angle $\alpha$ in its invariant plane combined with a reflection in that plane. Hence the name {\em rotoreflection} for such an isometry.
\\
\\
In this section we prove the {\sc Euler--Rodrigues} formula for 3D rotations. This formula arises from the general 4D rotation formula by setting $p = a$, $q = -b$, $r = -c$, $s = -d$ in Eq. \ref{Label2}.
Expressing this in quaternion terms one obtains the {\sc Hamilton--Cayley} formula $P' = QPQ^{-1}$, or expanded: $u'+x'i+y'j+z'k = (a+bi+cj+dk)(u+xi+yj+zk)(a-bi-cj-dk)$.
\\
\\
{\sc Proof:} The general 3D rotation matrix can be obtained from the general 4D rotation matrix by putting
$a_{00} = 1$, $a_{10} = a_{20} = a_{30} = 0$, $a_{01} = a_{02} = a_{03} = 0$. Doing so one obtains
\\
\begin{equation}
\label{Label6a}
A =
\left[
\begin{array}{rrrr}
                         1  &        0  &        0  &        0  \\
\noalign{\medskip}       0  &   a_{11}  &   a_{12}  &   a_{13}  \\
\noalign{\medskip}       0  &   a_{21}  &   a_{22}  &   a_{23}  \\
\noalign{\medskip}       0  &   a_{31}  &   a_{32}  &   a_{33}
\end{array}
\right].
\end{equation}
Its associate matrix reads
\\
\begin{equation}
\label{Label7a}
M = \frac{1}{4}
\left[
\begin{array}{llll}
                    1+a_{11}+a_{22}+a_{33} & +a_{23}-a_{32} & +a_{31}-a_{13} & +a_{12}-a_{21} \\
\noalign{\medskip}  a_{32}-a_{23} & -1-a_{11}+a_{22}+a_{33} & -a_{21}-a_{12} & -a_{31}-a_{13} \\
\noalign{\medskip}  a_{13}-a_{31} & -a_{21}-a_{12} & -1+a_{11}-a_{22}+a_{33} & -a_{32}-a_{23} \\
\noalign{\medskip}  a_{21}-a_{12} & -a_{31}-a_{13} & -a_{32}-a_{23} & -1+a_{11}+a_{22}-a_{33}
\end{array}
\right]
\end{equation}
\\
Comparing the off--diagonal elements with the corresponding elements in Eq. \ref{Label4} one infers
\\
\begin{equation}
\label{Label8a}
aq = -bp,\ \ ar = -cp, \ \ as = -dp,
\end{equation}
\begin{equation}
\label{Label9a}
cq = br, \ \ dq = bs, \ \ dr = cs.	
\end{equation}
Let $\alpha$ be the rotation angle of the 3D rotation. The trace of its matrix is $a_{11} + a_{22} + a_{33} = 2\cos \alpha + 1$, and the diagonal elements of the assoicate matrix are
\\
\begin{equation}
\label{Label10a}
ap = (a_{11} + a_{22} + a_{33} + 1)/4 \ \ = \ \ (2\cos \alpha + 2)/4	
\end{equation}
\begin{eqnarray}
\label{Label11a}
bq = (a_{11} + a_{22} + a_{33} - 1 - 2a_{11})/4 \ & = \ (2\cos \alpha - 2a_{11})/4  \\ 
\label{Label12a}
cr = (a_{11} + a_{22} + a_{33} - 1 - 2a_{22})/4 \ & = \ (2\cos \alpha - 2a_{22})/4  \\
\label{Label13a}
ds = (a_{11} + a_{22} + a_{33} - 1 - 2a_{33})/4 \ & = \ (2\cos \alpha - 2a_{33})/4
\end{eqnarray}
\\
From Eq. \ref{Label10a} one deduces
\begin{equation}
\label{Label14a}
ap \geq 0.
\end{equation}
The quantities $a_{11}$, $a_{22}$, $a_{33}$ are the cosines of the angles through which the $+OX$, $+OY$ and $+OZ$ half--axes are displaced. These angles are at most $\alpha$, therefore their cosines are at least $\cos \alpha$. For a proof of this see Section \ref{SectionProofs}. We conclude that 
\begin{equation}
\label{Label15a}
bq \leq 0, \ \ cr \leq 0, \ \ ds \leq 0.
\end{equation}
We recall the conditions
\\
\begin{equation}
\label{Label16a}	
a^2 + b^2 + c^2 + d^2 = 1,\ \ p^2 + q^2 + r^2 + s^2 = 1.
\end{equation}
From Eq. \ref{Label8a} follows the equality
\begin{equation}
\label{Label17a}	
a^2(q^2 + r^2 + s^2) = p^2(b^2 + c^2 + d^2),
\end{equation}
which is because of Eq. \ref{Label16a} the same as $a^2(1 - p^2) = p^2(1 - a^2)$, or equivalently: $a^2 = p^2$.
\\
Taking into account Eq. \ref{Label14a} we conclude $a = p$.
\\
Substitute $a$ for $p$ in Eq. \ref{Label8a} to obtain $aq = -ab$, $ar = -ac$, $as = -ad$.
\\
If $a \neq 0$ then $q = -b$, $r = -c$, $s = -d$.
\\
If $a = 0$ then also $p = 0$. In that case one has
\begin{equation}
\label{Label18a}
b^2 + c^2 + d^2 = 1,\ \ q^2 + r^2 + s^2 = 1.
\end{equation}
Eq. \ref{Label9a} is equivalent to $(b, c, d) \propto (q, r, s)$.
Combining this with Eq. \ref{Label18a} and Eq. \ref{Label15a} we conclude that $q = -b$, 
$r = -c$, $s = -d$ also for $a = 0$. End of proof.
\\

\section{The Euler--Rodrigues formula}
\label{Section4}

By substituting $a$, $-b$, $-c$, $-d$ for $p$, $q$, $r$, $s$ in Eq. \ref{Label2} and deleting the leftmost column and uppermost row one obtains the {\sc Euler--Rodrigues} formula in its usual form:

\begin{equation}
\label{Label19a}
\left[
\begin{array}{rrr}
                    a^2+b^2-c^2-d^2 \ \ & -2ad+2bc         \ \ &  2ac+2bd \\
\noalign{\medskip}  2ad+2bc         \ \ &  a^2-b^2+c^2-d^2 \ \ & -2ab+2cd \\
\noalign{\medskip} -2ac+2bd         \ \ &  2ab+2cd         \ \ &  a^2-b^2-c^2+d^2
\end{array}
\right]
\end{equation}
with the relation $a^2 + b^2 + c^2 + d^2 = 1$.
\\
\\
This formula was first discovered by {\sc Euler} ([EULE~1770]) and later rediscovered independently by {\sc Rodrigues} ([RODR~1840]), hence the terms {\em {\sc Euler} parameters} for the parameters $a,\ b,\ c,\ d$ and {\em {\sc Euler--Rodrigues} formula} (or {\em formulae}) for the rotation matrix in terms of $a,\ b,\ c,\ d$.
\\

\subsection{Determination of the Euler parameters for a given 3D rotation}

By substituting $a$, $-b$, $-c$, $-d \ $ for $p$, $q$, $r$, $s \ $ in the associate matrix (Eq. \ref{Label5}), or directly from the {\sc Euler--Rodrigues} formula and the accompanying relation, one obtains the ten equations below in the four unknowns $a$, $b$, $c$, $d$:

\begin{eqnarray}
\label{Label20}
a^2 & = & (1 + a_{11} + a_{22} + a_{33})/4,\\	
\label{Label21}
b^2 & = & (1 + a_{11} - a_{22} - a_{33})/4,\\	
\label{Label22}
c^2 & = & (1 - a_{11} + a_{22} - a_{33})/4,\\	
\label{Label23}
d^2 & = & (1 - a_{11} - a_{22} + a_{33})/4,
\end{eqnarray}
\begin{eqnarray}
\label{Label24}
ab & = & (a_{32} - a_{23})/4,\\
\label{Label25}
ac & = & (a_{13} - a_{31})/4,\\
\label{Label26}
ad & = & (a_{21} - a_{12})/4,
\end{eqnarray}
\begin{eqnarray}
\label{Label27}
cd & = & (a_{32} + a_{23})/4,\\
\label{Label28}
db & = & (a_{13} + a_{31})/4,\\
\label{Label29}
bc & = & (a_{21} + a_{12})/4,
\end{eqnarray}
which are consistent and admit a pair of real solutions if and only if the given matrix is indeed a 3D rotation matrix. 
\\
\\
Assuming that the given matrix $a_{11} \ldots a_{33}$ is indeed a 3D rotation matrix, one first determines the absolute values of $a$, $b$, $c$, $d$ from Eqs. \ref{Label20} $\ldots$  \ref{Label23}; next one arbitrarily fixes the sign of one of the nonzero quantities $a$, $b$, $c$, $d$ (at least one of them is nonzero); finally one determines the signs of the other ones from the remaining six equations. 
\\
There always exists a subset of Eqs. \ref{Label20} $\ldots$ \ref{Label29} from which one can determine $a$, $b$, $c$, $d$ save for their signs; the remaining equations provide opportunities for cross--checking the calculations and for averaging numerical errors. Think in this connection of ill--conditioned data and of input rotation matrices affected with errors.
\\
In practice one uses Eqs. \ref{Label24} $\ldots$ \ref{Label26} for rotations through small angles, and Eqs. \ref{Label27} $\ldots$ \ref{Label29} for rotations through angles close to $\pi$.
\\
These procedures are well--known in the fields of aviation and astrodynamics; see for instance [BATT~1999] and [SHEP~1978].

\section{3D Rotoreflections}
\label{Section3}

In this section we prove the counterpart of the {\sc Euler--Rodrigues} formula for 3D rotoreflections. This formula arises from the general 4D rotation formula by setting $p = -a$, $q = b$, $r = c$, $s = d$ in Eq. \ref{Label2}.
Expressing this in quaternion terms one obtains the {\sc Hamilton--Cayley}--like formula $P' = -QPQ^{-1}$, or expanded: $u'+x'i+y'j+z'k = -(a+bi+cj+dk)(u+xi+yj+zk)(a-bi-cj-dk)$.
\\
\\
{\sc Proof:} The general 3D rotoreflection matrix can be obtained from the general 4D rotation matrix by putting
$a_{00} = -1$, $a_{10} = a_{20} = a_{30} = 0$, $a_{01} = a_{02} = a_{03} = 0$. Doing so one obtains
\\
\begin{equation}
\label{Label6b}
A =
\left[
\begin{array}{rrrr}
                        -1  &        0  &        0  &        0  \\
\noalign{\medskip}       0  &   a_{11}  &   a_{12}  &   a_{13}  \\
\noalign{\medskip}       0  &   a_{21}  &   a_{22}  &   a_{23}  \\
\noalign{\medskip}       0  &   a_{31}  &   a_{32}  &   a_{33}
\end{array}
\right].
\end{equation}
Its associate matrix reads
\\
\begin{equation}
\label{Label7b}
M = \frac{1}{4}
\left[
\begin{array}{llll}
                  -1+a_{11}+a_{22}+a_{33} & +a_{23}-a_{32} & +a_{31}-a_{13} & +a_{12}-a_{21} \\
\noalign{\medskip}  a_{32}-a_{23} & 1-a_{11}+a_{22}+a_{33} & -a_{21}-a_{12} & -a_{31}-a_{13} \\
\noalign{\medskip}  a_{13}-a_{31} & -a_{21}-a_{12} & 1+a_{11}-a_{22}+a_{33} & -a_{32}-a_{23} \\
\noalign{\medskip}  a_{21}-a_{12} & -a_{31}-a_{13} & -a_{32}-a_{23} & 1+a_{11}+a_{22}-a_{33}
\end{array}
\right]
\end{equation}
Comparing the off--diagonal elements with the corresponding elements in Eq. \ref{Label4} one infers
\\
\begin{equation}
\label{Label8b}
aq = -bp,\ \ ar = -cp, \ \ as = -dp,
\end{equation}
\begin{equation}
\label{Label9b}
cq = br, \ \ dq = bs, \ \ dr = cs.	
\end{equation}
Let $\alpha$ be the rotation angle of the 3D rotoreflection. The trace of its matrix is $a_{11} + a_{22} + a_{33} = 2\cos \alpha - 1$, and the diagonal elements of the assoicate matrix are
\\
\begin{equation}
\label{Label10b}
ap = (a_{11} + a_{22} + a_{33} - 1)/4 \ \ = \ \ (2\cos \alpha - 2)/4	
\end{equation}
\begin{eqnarray}
\label{Label11b}
bq = (a_{11} + a_{22} + a_{33} + 1 - 2a_{11})/4 \ & = \ (2\cos \alpha - 2a_{11})/4  \\ 
\label{Label12b}
cr = (a_{11} + a_{22} + a_{33} + 1 - 2a_{22})/4 \ & = \ (2\cos \alpha - 2a_{22})/4  \\
\label{Label13b}
ds = (a_{11} + a_{22} + a_{33} + 1 - 2a_{33})/4 \ & = \ (2\cos \alpha - 2a_{33})/4
\end{eqnarray}
\\
From Eq. \ref{Label10b} one deduces
\begin{equation}
\label{Label14b}
ap \leq 0.
\end{equation}
The quantities $a_{11}$, $a_{22}$, $a_{33}$ are the cosines of the angles through which the $+OX$, $+OY$ and $+OZ$ half--axes are displaced. These angles are at least $\alpha$, therefore their cosines are at most $\cos \alpha$. For a proof of this see Section \ref{SectionProofs}. We conclude that 
\begin{equation}
\label{Label15b}
bq \geq 0, \ \ cr \geq 0, \ \ ds \geq 0.
\end{equation}
We recall the conditions
\\
\begin{equation}
\label{Label16b}	
a^2 + b^2 + c^2 + d^2 = 1,\ \ p^2 + q^2 + r^2 + s^2 = 1.
\end{equation}
From Eq. \ref{Label8b} follows the equality
\begin{equation}
\label{Label17b}	
a^2(q^2 + r^2 + s^2) = p^2(b^2 + c^2 + d^2),
\end{equation}
which is because of Eq. \ref{Label16b} the same as $a^2(1 - p^2) = p^2(1 - a^2)$, or equivalently: $a^2 = p^2$.
\\
Taking into account Eq. \ref{Label14b} we conclude $a = -p$.
\\
Substitute $-a$ for $p$ in Eq. \ref{Label8b} to obtain $aq = ab$, $ar = ac$, $as = ad$.
\\
If $a \neq 0$ then $q = b$, $r = c$, $s = d$.
\\
If $a = 0$ then also $p = 0$. In that case one has
\begin{equation}
\label{Label18b}
b^2 + c^2 + d^2 = 1,\ \ q^2 + r^2 + s^2 = 1.
\end{equation}
Eq. \ref{Label9b} is equivalent to $(b, c, d) \propto (q, r, s)$.
Combining this with Eq. \ref{Label18b} and Eq. \ref{Label15b} we conclude that $q = b$, 
$r = c$, $s = d$ also for $a = 0$. End of proof.
\\
\\
The counterpart of the {\sc Euler--Rodrigues} formula for 3D rotoreflections reads
\begin{equation}
\label{Label19b}
\left[
\begin{array}{rrr}
                   -a^2-b^2+c^2+d^2 \ \ & +2ad-2bc         \ \ & -2ac-2bd \\
\noalign{\medskip} -2ad-2bc         \ \ & -a^2+b^2-c^2+d^2 \ \ & +2ab-2cd \\
\noalign{\medskip} +2ac-2bd         \ \ & -2ab-2cd         \ \ & -a^2+b^2+c^2-d^2
\end{array}
\right]
\end{equation}
with the relation $a^2 + b^2 + c^2 + d^2 = 1$.
\\
\\

\section{Proofs of the angle inequalities for 3D rotations and rotoreflections}
\label{SectionProofs}

The matrix of a 3D rotation around the $Z$ axis reads
\[
R =
\left[
\begin{array}{rrr}
                    \cos\alpha  & -\sin\alpha  &  0  \\
\noalign{\medskip}  \sin\alpha  &  \cos\alpha  &  0  \\
\noalign{\medskip}           0  &           0  &  1
\end{array}
\right]
\]

Let $\beta$ be the angle through which the ray from $O$ through a point $(x, y, z)$ is displaced. An easy calculation yields (putting $\rho^2 = x^2 + y^2$)
\[
\cos\beta = (\rho^2\cos\alpha + z^2)/(\rho^2 + z^2),
\]
\[
\cos\beta - \cos\alpha = (\rho^2\cos\alpha + z^2 - \rho^2\cos\alpha - z^2\cos\alpha)/(\rho^2 + z^2) = z^2(1-\cos\alpha)/(\rho^2 + z^2) \geq 0.
\]
\\

The matrix of a 3D rotoreflection with the $Z$ axis as its axis reads

\[
R =
\left[
\begin{array}{rrr}
                    \cos\alpha  & -\sin\alpha  &  0  \\
\noalign{\medskip}  \sin\alpha  &  \cos\alpha  &  0  \\
\noalign{\medskip}           0  &           0  & -1
\end{array}
\right]
\]

Let $\beta$ be the angle through which the ray from $O$ through a point $(x, y, z)$ is displaced. An easy calculation yields (putting $\rho^2 = x^2 + y^2$)
\[
\cos\beta = (\rho^2\cos\alpha - z^2)/(\rho^2 + z^2),
\]
\[
\cos\beta - \cos\alpha = (\rho^2\cos\alpha - z^2 - \rho^2\cos\alpha - z^2\cos\alpha)/(\rho^2 + z^2) = -z^2(1+\cos\alpha)/(\rho^2 + z^2) \leq 0.
\]
\\
\newpage
\section*{Literature}

\begin{tabbing}
WWWJ\=WWWJ \=DUMMY             \kill
BATT \>1999 \>{\sc Richard H. Battin:} \\
            \>\>An introduction to the mathematics and methods of astrodynamics, revised edition.\\  
            \>\>AIAA Education Series\\  
            \>\>{\em American Institute of Aeronautics and Astronautics, Inc.}, 1999\\  
            \>\>\\  
EULE \>1770 \>{\sc L.\ Euler}: \\
            \>\>Problema algebraicum ob affectiones prorsus singulares memorabile. \\
            \>\>{\em Commentatio 407 indicis Enestr\oe miani, Novi commentarii} \\
            \>\>{\em academi\ae \ scientiarum Petropolitan\ae}, 15, (1770), 1771, \\
            \>\>p.75-106, reprinted in L.Euleri Opera Omnia, 1$^{st}$ series, \\
            \>\>Vol. 6, p.287-315 \\
            \>\>\\
MEBI \>2005 \>{\sc Johan Ernest Mebius}:\\
            \>\>A matrix--based proof of the quaternion representation theorem for four--dimensional rotations.\\
            \>\>{\em arXiv General Mathematics} 2005 \\
            \>\>{\tt http://www.arxiv.org/PS\_cache/math/pdf/0501/0501249.pdf}\\
            \>\>\\
RODR \>1840 \>{\sc Olinde Rodrigues}: \\
            \>\>Des lois g\a'{e}om\a'{e}triques qui r\a'{e}gissent les d\a'{e}placements d'un syst\a`{e}me solide dans l'espace, \\
            \>\>et de la variation des coordonn\a'{e}es provenant de ces d\a'{e}placements consid\a'{e}r\a'{e}s \\
            \>\>ind\a'{e}pendamment des causes qui peuvent les produire.\\
            \>\>{\em Journal de Math\a'{e}matiques} 5, 1840, 380-440\\
            \>\>\\  
SHEP \>1978 \>{\sc Stanley W. Sheppard:} \\
            \>\>Quaternion from Rotation Matrix.\\
            \>\>{\em Journal of Guidance and Control}, Vol. 1 nr 3, May--June 1978, pp. 223-224 \\
            \>\>\\  
\end{tabbing}

\end{document}